\documentclass{amsart}
\usepackage{amssymb}
\usepackage{cite}

\newtheorem{theorem}{Theorem}[section]

\theoremstyle{definition}

\theoremstyle{remark}
\newtheorem{remark}[theorem]{Remark}

\numberwithin{equation}{section}

\newcommand*{\eqb}{\begin{eqnarray}}
\newcommand*{\eqe}{\end{eqnarray}}
\newcommand*{\al}{\alpha}
\newcommand*{\N}{\mathbb{N}}
\newcommand*{\C}{\mathbf{C}}
\newcommand*{\R}{\mathbb{R}}
\newcommand*{\lap}{\mathcal{L}}
\newcommand*{\eb}{\begin{equation}}
\newcommand*{\ee}{\end{equation}}
\newcommand*{\ebnn}{\begin{equation*}}
\newcommand*{\eenn}{\end{equation*}}
\newcommand*{\abs}[1]{\left| #1 \right|}
\newcommand*{\ind}{\mathbf{1}}

\newcommand*{\ro}{\varrho}
\newcommand*{\ex}{\mathbb{E}}

\begin{document}
\title[Stochastic representation of fractional subdiffusion equation]{Stochastic
representation of fractional subdiffusion equation.
The case of infinitely divisible waiting times, L\'evy noise and space-time-dependent coefficients}
\author[Marcin Magdziarz]{Marcin Magdziarz
\and Tomasz Zorawik}
\address{Hugo Steinhaus Center,
Institute of Mathematics and Computer Science,
Wroclaw University of Technology,
50-370 Wroclaw, Poland}
\email{marcin.magdziarz@pwr.wroc.pl}

\subjclass[2010]{\subjclass[2000]{Primary 60G22; Secondary 60G51}}

\begin{abstract}
In this paper we analyze fractional Fokker-Planck equation describing
subdiffusion in the general infinitely divisible (ID) setting.
We show that in the case of space-time-dependent drift and diffusion and time-dependent jump coefficient,
the corresponding stochastic process can be obtained by subordinating two-dimensional system of
Langevin equations driven by appropriate Brownian and L\'evy noises.
Our result solves the problem of stochastic representation of subdiffusive Fokker-Planck dynamics
in full generality.
\end{abstract}

\maketitle
\section{Preliminaries}
Subdiffusion processes are characterized by the asymptotic power-law behavior of the variance $\mathrm{Var}(X(t))\sim ct^\alpha$ as $t\rightarrow \infty$, where $0<\alpha<1$.
Physical equation describing
subdiffusion in the presence of space-dependent
force $F(x)$ is the fractional Fokker-Planck equation (FFPE) \cite{Metzler, MBK,MBK2, metzler2000}
\eqb
\label{ffp_space}
\frac{\partial w(x,t)}{\partial t}=
{_0D_t^{1-\al}}\left[ -\frac{\partial}{\partial x} F(x)
+\frac{\sigma^2}{2} \frac{\partial^2}{\partial x^2}  \right]w(x,t),
\eqe
with the scale parameter $\sigma>0$ and the initial condition $w(x,0)=\delta(x)$.
Here
\[
_0D_t^{1-\al} f(t)=\frac{1}{\Gamma(\al)} \frac{d}{dt}\int_0 ^t (t-s)^{\al-1} f(s)ds,
\]
$0<\al<1, \;f\in C^1([0,\infty)),$ is the fractional derivative of the Riemann-Liouville type \cite{samko}.
Eq. \eqref{ffp_space} specifies probability density function (PDF) $w(x,t)$
of subdiffusion process.
It was derived in \cite{MBK} in the framework of continuous-time random walk
with heavy-tailed waiting times \cite{Meer6,Meer_Book}.

Stochastic process corresponding to \eqref{ffp_space}
has the subordination form
\[
Y(t)=X(S_\al(t)), \;\;\;t\geq 0,
\]
where $X$ is given by the following stochastic differential equation
\eqb
\label{SDE}
dX(t)={F(X(t))}dt+\sigma dB(t),\;\; X(0)=0,
\eqe
driven by Brownian motion $B$.
Moreover, $S_{\al}$ is the inverse $\al$-stable subordinator
\eqb
\label{inv_sub}
S_{\al}(t)=\inf\{\tau>0: T_{\al}(\tau)>t \},
\eqe
which is assumed independent of $B$. Here, $T_{\al}$ is the $\al$-stable
subordinator \cite{Sato,Janicki_Weron} with Laplace transform
$\mathbb{E}\left[ e^{-uT_{\al}(\tau)}   \right]=e^{-t u^\al}$,
$0<\al<1$. As shown in \cite{MWW,Meer1}, PDF of $Y(t)=X(S_\al(t))$ solves FFPE \eqref{ffp_space}. The concept of subordination in terms of coupled Langevin equation was
originally introduced in \cite{fogedby}.

FFPE describing subdiffusion in the presence of time-dependent force $F(t)$
was derived in \cite{Klaf_Sok} and has the form
\eqb
\label{ffp_time}
\frac{\partial w(x,t)}{\partial t}=
\left[ -F(t)\frac{\partial}{\partial x}
+\frac{\sigma^2}{2} \frac{\partial^2}{\partial x^2}  \right]{_0D_t^{1-\al}}w(x,t),
\eqe
$w(x,0)=\delta(x)$.
Compared to (\ref{ffp_space}), we can observe that the fractional operator ${_0D_t^{1-\al}}$ in the above
equation does not act on the force $F(t)$, so the force is not influenced by the change of time.

Stochastic process corresponding to \eqref{ffp_time} has the form \cite{Magdziarz1}
\eqb
\label{stoch_repr_time}
Y(t)=\int_0^t F(u)dS_{\al}(u)+\sigma B(S_{\al}(t)),\;\;t\geq 0,
\eqe
meaning that the PDF of $Y(t)$ solves \eqref{ffp_time}.
$Y(t)$ in \eqref{stoch_repr_time} can be equivalently represented in the form of subordination
\eqb
\label{stoch_repr_time_2}
Y(t)=X(S_{\al}(t)),
\eqe
where $X$ is given by
\[
d{X}(t)=F(T_{\al}(t))dt+\sigma dB(t),\;\; {X}(0)=0,
\]
and $S_\al$ is the inverse subordinator independent of $B$. Note that in this case
processes $X$ and $S_\al$ are not independent, since $T_\al$ and its inverse $S_\al$
are strongly dependent, see formula \eqref{inv_sub}.

Further relevant studies of FFPE in the $\al$-stable case include:
fractional Cauchy problems \cite{Meer2,Chen,Mijena}, path properties
of fractional diffusion \cite{Meer3,Nane,Magdziarz2}, Pearson diffusion \cite{Leonenko}
and the case of general forces \cite{Straka,Chinczycy}.
In the recent paper \cite{MGZ} the case of a general space-time dependent force and diffusion coefficients is considered.

Extension of \eqref{stoch_repr_time} (or equivalently \eqref{stoch_repr_time_2})
to the ID case is the following:
by $T_\Psi(t)$, $t\geq 0$ we denote a subordinator (strictly increasing L\'evy process)
with Laplace transform \cite{Bertoin}
\[
\mathbb{E}\left( e^{-uT_\Psi(t)}   \right)=e^{-t \Psi(u)},
\]
Here $\Psi(u)$ is the Laplace exponent given by
\[
\Psi(u)=\lambda u+\int_0^\infty (1-e^{-ux})\nu (dx).
\]
We assume for simplicity that $\lambda=0$.
The measure $\nu(dx)$ is the L\'evy measure supported in $(0,\infty)$
satisfying $\int_{(0,\infty)}\min\{1,x\}\nu(dx)<\infty$.
To exclude the case of compound Poisson processes, we
will further assume that $\nu((0,\infty))=\infty$.
The corresponding first-passage-time process
\eqb
\label{inv_subordinator}
S_\Psi(t)=\inf\{\tau>0: T_\Psi(\tau)>t \},\;\;\;\;\; t\geq 0.
\eqe
is called {\em inverse subordinator.}
Inverse subordinators play an important role
in probability theory \cite{Bertoin,Bertoin2},
finance \cite{Winkel} and physics \cite{Jurlewicz_Weron,SWW}.
Since $\nu((0,\infty))=\infty$, the trajectories of $S_\Psi(t)$ are continuous,
which is relevant in physical applications.

Now, to extend \eqref{stoch_repr_time} and \eqref{stoch_repr_time_2} to the ID
case, we replace the inverse stable subordinator $S_\al$ with $S_\Psi$
and obtain
\eqb
\label{stoch_repr_time_ID}
Y(t)=\int_0^t F(u)dS_{\Psi}(u)+\sigma B(S_{\Psi}(t)),\;\;t\geq 0,
\eqe
or equivalently
\eqb
\label{stoch_repr_time_2_ID}
Y(t)=X(S_{\Psi}(t)),
\eqe
In this case $X$ is given by
\eb
\label{id_x_sde}
d{X}(t)=F(T_{\Psi}(t))dt+\sigma dB(t),\;\; {X}(0)=0,
\ee
where the subordinator $T_\Psi$ is independent of Brownian motion $B$.
The process $Y$ defined in \eqref{stoch_repr_time_ID} and \eqref{stoch_repr_time_2_ID}
describes subdiffusion with ID waiting times in the presence
of time-dependent force $F(t)$ \cite{Klaf_Sok,Magdziarz3}.
The corresponding fractional Fokker-Planck equation has the form \cite{Klaf_Sok,Magdziarz3}
\eqb
\label{ffp_id}
\frac{\partial w(x,t)}{\partial t}=
\left[ -F(t)\frac{\partial}{\partial x}
+\frac{\sigma^2}{2} \frac{\partial^2}{\partial x^2}  \right] {\Phi}_t w(x,t).
\eqe
The integro-differential operator ${\Phi}_t$ is defined as
\eqb
\label{operator}
\Phi_t f(t)=\frac{d}{dt} \int_0^t M(t-y)f(y)dy,
\eqe
where $f$ is a sufficiently smooth function.
Moreover, the kernel $M(t)$ is defined via its
Laplace transform
\eqb
\label{kernel}
\lap[M](u)=\int_0^\infty e^{-ut}M(t)dt=\frac{1}{\Psi(u)},
\eqe
where $\Psi(u)$ is the Laplace exponent of the underlying ID distribution of waiting times.

Except from the previously mentioned $\al$-stable case, other important in applications examples include:
tempered-stable case \cite{Rosinski,Magdziarz3,Meer7,Meer2a}, which corresponds to
$\Psi(u)=(u+\lambda)^\alpha-\lambda^\alpha$ with $\lambda>0$ and $0<\alpha<1$;
distributed order case \cite{Kochubei,Chechkin,Meer9},
which corresponds to $\Psi(u)=\int_0^\infty (1-e^{-ux})\nu (dx)$ with
$\nu(t,\infty)=\int_0^1 t^{-\beta} \mu(d\beta)$ (here $\beta\in (0,1)$ and $\mu$ is
some distribution supported in $[0,1]$).

Another extension of \eqref{stoch_repr_time_2_ID} is possible when we add a L\'evy noise $dL(t)$ into Eq. \eqref{id_x_sde}. For the L\'evy process $L(t)$ we have the L\'evy-Khintchine formula
\eb
\label{lc}
\nonumber
\ex e^{iuX(t)}=\exp\left[t\left(ibu-\frac{a^2u^2}{2}+ \int_{\R\setminus \{0\}}\left(e^{iuy}-1-iuy\ind_{\{|y|<1\}}\right)\mu(dy)\right)\right],
\ee
where $b$ is a drift parameter, $a^2$ is connected with Brownian part of $L(t)$ and $\mu$ is the L\'evy measure with $\mu(dy)$ being the intensity of jumps of size $y$. Since in  our stochastic differential equations the drift and Brownian part already exist, we may assume without loss of generality that $b=a=0$.

The idea of stochastic representation of FFPE has become very popular in a physical society. In a very recent paper \cite{Zhang} authors find the stochastic representation in the case of tempered - stable subordination with space-time dependent coefficients.
In the next section we significantly extend the results from \cite{Magdziarz3,MGZ, Zhang}: we consider the general case of space-time-dependent drift and diffusion coefficients, ID waiting times and L\'evy noise. Moreover, we analyze the Fokker-Planck equation resulting from adding the L\'evy noise $E(t)dL(t)$ to Eq. \eqref{id_x_sde}.
The general idea of the proof of the main theorem is similar to \cite{MGZ}, however we focus our attention on necessary changes coming from the ID subordinator and the L\'evy noise. The most important differences are underlined in the proof of the main theorem. We also put our effort into delivering all necessary mathematical details that were lacking in the physical paper \cite{MGZ}.
\section{Main result}
Now, let us look at the following fractional Fokker-Planck equation
\eb
\label{ffp_general}
\begin{split}
\frac{\partial}{\partial t} &w(x,t)=\left[ -\frac{\partial}{\partial x} F(x,t)
+\frac{1}{2} \frac{\partial^2}{\partial x^2}\sigma^2(x,t)  \right]{\Phi_t} w(x,t)\\
&+
\int_{\R\setminus\{0\}}\left( \Phi_t w(r,t)\bigg|_{r=x-E(t)y}
-\Phi_t w(x,t)
+E(t)y\frac{\partial}{\partial x}\Phi_t w(x,t)\ind_{\{|y|<1\}}\right)\nu(dy),
\end{split}
\ee
$w(x,0)=\delta(x)$. It describes the temporal evolution of $w(x,t)$ -- PDF of
some anomalous diffusion process with ID waiting times, which is subjected to space-time-dependent force $F(x,t)$ and diffusion $\sigma(x,t)$ and with time-dependent jump coefficient $E(t)$.
In the theorem below, which is the main result of the paper, we derive the stochastic process corresponding to \eqref{ffp_general}.
The result encompasses the whole family of ID waiting times,
arbitrary space-time-dependent drift and diffusion coefficients and L\'evy jumps with time-dependent jump coefficient. From now on for the left limit process we will use the notation $X^-(t)=\lim_{s\rightarrow t^-}X(s)$.
\begin{theorem}
Let $T_\Psi (\tau)$ be the subordinator with the Laplace exponent $\Psi(u)$ and $S_\Psi (t)$ its inverse.
Assume that the standard Brownian motion $B(t)$, the L\'evy process $L(t)$  with L\'evy measure $\mu$ and $T_\Psi(\tau)$ are independent. Let $Y(t)$ be the solution of the stochastic equation
\eqb
\label{s_representation}
&&dY(t)=F(Y^-(t),T^-_\Psi(t))dt+\sigma(Y^-(t),T^-_\Psi(t))dB(t)+E(T^-_\Psi(t))dL(t),\; t\geq 0,\\
&&Y(0)=0,\;\;\;T_\Psi(0)=0, \nonumber
\eqe
where the functions $F(x,t)$, $\sigma(x,t)$, $E \in \C^2(\R^2)$ satisfy the Lipschitz condition. We assume that the PDF of the process $(Y(t),T_\Psi(t))$ - $p_t(y,z)$ exists and so do $\frac{\partial}{\partial t} p_t(y,z)$, $\frac{\partial}{\partial y} p_t(y,z)$ and $\frac{\partial^2}{\partial y^2} p_t(y,z)$. Additionally we require that:
\eb
\label{dt_integrable}
\int_0^{t_0}\int_0^\infty\abs{\frac{\partial}{\partial t} p_s(x,t)}dsdt<\infty
\ee
for each $t_0>0$,
\eb
\label{dx_integrable}
\int_{x_1}^{x_2}\int_0^\infty\abs{\frac{\partial}{\partial x}p_s(x,t)}dsdx<\infty, \quad \int_{x_1}^{x_2}\int_0^\infty\abs{\frac{\partial^2}{\partial x^2}p_s(x,t)}dsdx<\infty
\ee
for each $x_1$, $x_2\in \R$ and finally
\eb
\label{ds_integrable}
\int_0^\infty\int_{\R\setminus\{0\}}\abs{ p_s(x-E(t)y,t)
-p_s(x,t)
+\frac{\partial}{\partial x}(E(t)y) p_s(x,t)\ind_{\{|y|<1\}}(y)}\nu(dy)ds<\infty
\ee
for each $x\in \R$, $t>0$.  Then the PDF of the process $X(t)=Y^-(S_\Psi (t))$ is the weak solution of FFPE \eqref{ffp_general}, that is \eqref{ffp_general} holds pointwise for $t>0$ with the required initial distribution for $t=0$.
\end{theorem}
\begin{remark}
If we impose additional conditions on the coefficients $F,\sigma, E$,  the L\'evy measures $\mu$ and $\nu$ (of $L(t)$ and $T_\Psi(t)$, respectively) namely: $F,\sigma \in C^3_b(\R^2)$, $E \in C^3_b(\R)$  with bounded partial derivatives of the order $1-3$, $D(y,z)\neq 0$ for all $y,z\in \R$  and $\int_0^\infty x^p \nu(dx)<\infty$, $\int_{\R} \abs{x}^p \mu(dx)<\infty$ for all $p\geq 2$, then the density $p_t(y,z)$  exists, see Theorem 2-14 in \cite{Bichteler}. Other, less strict non-degeneracy conditions are also discussed there. The regularity assumptions of $p_t(y,z)$ are easy to check when $F(x,t)$ and $\sigma(x,t)$ depend only on $x$ and the processes $Y(t)$ and $T_\Psi(t)$ are independent. Notice that we do not assume anything about the PDF of $X(t)$.
\end{remark}
\begin{proof}
In the proof we will essentially extend the techniques used in  \cite{MGZ,Straka,Magdziarz1}.
Eq. \eqref{s_representation} is equivalent to the following system of stochastic equations
\eb
\begin{split}
\label{SDE_matrix}
\left(\begin{array}{c}
dY(t)\cr
dZ(t)\cr
\end{array}\right)&=\left(\begin{array}{c}
F(Y^-(t),Z^-(t))\cr
0\cr
\end{array}\right)dt+
\left(\begin{array}{c}
\sigma(Y^-(t),Z^-(t))\cr
0\cr
\end{array}\right)dB(t)\\
&+\left(\begin{array}{c}
E(Z^-(t))dL(t)\cr
dT_\Psi(t)\cr
\end{array}\right).
\end{split}
\ee
This system is subjected to Brownian and L\'evy noise, therefore the infinitesimal generator for the process $\left(Y(t),Z(t)\right)$
operates on functions $f\in \C_0^2(\R^2)$ in the following way (see Theorem 6.7.4 in \cite{Applebaum}):
\eb
\begin{split}
\label{generator}
Af(y,z)&= F(y,z)\frac{\partial}{\partial y} f(y,z)+ \frac{1}{2}\sigma^2(y,z)\frac{\partial^2}{\partial y^2} f(y,z)\\ &+\int_{\R\setminus\{0\}}\left(f(y+E(z)x,z)-f(y,z)-E(z)x\frac{\partial}{\partial y}f(y,z)\ind_{\{\abs{x}<1\}}\right)\mu(dx)\\
&+\int_0^\infty \left[f(y,z+u)-f(y,z) \right]\nu(du),
\end{split}
\ee
where $\nu$ is the L\'evy measure of $Z(t)$ (or equivalently $T_\Psi(t)$). We have two separate integrals here, because $L(t)$ and $T_\Psi(t)$ are assumed to be independent which implies that $\ro$ - L\'evy measure of the two-dimensional  process $(L(t),T_\Psi(t))$ is supported
on the axes of $\R^{2}$, that is $\ro(A)=\mu(\{(y,z)\in \R^2: (y,0)\in  A\})+\nu(\{z\in \R:(0,z)\in A\})$ for $A\in Bor(\R^{2})$.
In the first step of the proof we will obtain a FFPE describing the temporal evolution of the density $p_t(y,z)$.
To do this we take advantage of the fact that
\eb
\label{q_adjoint}
\frac{\partial}{\partial t}p_t(y,z)=A^+p_t(y,z).
\ee
Here $A^+$ is $L^2$ Hermitian adjoint of $A$, meaning that it satisfies
the following relation for all test functions $f \in \C_c^\infty(\R^2)$:
\eb
\label{adjoint_relation}
\begin{split}
&\int_{\R^2} Af(y,z)p_t(y,z)dydz = \int_{\R^2}A^+p_t(y,z)f(y,z)dydz.
\end{split}
\ee
If we substitute \eqref{generator} into the above equation, then the
left-hand side consists of four summands. The first two summands are dealt with in \cite{MGZ}. Now we turn our attention to the third one - connected with an integration with respect to the measure $\mu$. Since $f$ has compact support we can apply Fubini's theorem to get
\eb
\begin{split}
\label{generator_jumps}
&\int_{\R^{2}}\int_{\R\setminus\{0\}}\\
&\quad\quad\quad\left(f(y+E(z)x,z)-f(y,z)-E(z)x\frac{\partial}{\partial y}f(y,z)\ind_{\{\abs{x}<1\}}\right)\mu(dx)p_t(y,z)d(y,z)\\
&=\int_{\R\setminus\{0\}}\int_{\R}\int_{\R}\\
&\quad\quad\quad\left(f(y+E(z)x,z)-f(y,z)-E(z)x)\frac{\partial}{\partial y}f(y,z)\ind_{\{\abs{x}<1\}}\right)p_t(y,z)dy dz\mu(dx).
\end{split}
\ee
For the inner integral (for each fixed $z$ and $x$) after a  substitution  and using integration by parts we get
\eb
\label{transformed2_generator_jumps}
\begin{split}
&\int_{\R}\left(f(y+E(z)x,z)-f(y,z)-E(z)x\frac{\partial}{\partial y}f(y,z)\ind_{\{\abs{x}<1\}}\right)p_t(y,z)dy \\
&=\int_{\R}f(y,z)\left(p_t(y-E(z)x,z)
-p_t(y,z)
+E(x)x\frac{\partial}{\partial y}p_t(y,z)\ind_{\{\abs{x}<1\}}\right)dy. \end{split}
\ee
 Substituting \eqref{transformed2_generator_jumps} back into \eqref{generator_jumps} and using Fubini's theorem again we obtain
\eb
\begin{split}
\label{generator_jumps_end}
&\int_{\R^{2}}\int_{\R\setminus\{0\}}\\
&\quad\quad\quad\left(f(y+E(z)x,z)-f(y,z)-E(z)x)\frac{\partial}{\partial y}f(y,z)\ind_{\{\abs{x}<1\}}\right)\mu(dx)p_t(y,z)d(y,z)\\
&=\int_{\R^{2}}f(y,z)\\
&\quad\quad\quad\int_{\R\setminus\{0\}}\left(p_t(y-E(z)x,z)
-p_t(y,z)
+E(z)x\frac{\partial}{\partial y}p_t(y,z)\ind_{\{\abs{x}<1\}}\right)\mu(dx)d(y,z).
\end{split}
\ee
Now let us put
\eb
G(u)=\nu ((u,\infty)).
\ee
Using the general formula for integration by parts (see (21.68) in  \cite{Hewitt}), we get for the last summand in Eq. \eqref{adjoint_relation}
\ebnn
\begin{split}
&\int_0^\infty \left[f(y,z+u)-f(y,z) \right] \nu(du) =  -\lim_{u \rightarrow \infty} \left(f(y,z+u)-f(y,z)\right)G(u)\\
&+\lim_{u \rightarrow 0} \left(f(y,z+u)-f(y,z)\right)G(u)
+\int_0^\infty \frac{\partial}{\partial u} f(y,z+u) G(u)du
\end{split}
\eenn
It is easy to see that the first limit vanishes. For the second limit we also have
\ebnn
\begin{split}
\lim_{u \rightarrow 0} \left(f(y,z+u)-f(y,z)\right)G(u)= \lim_{u \rightarrow 0} \frac{f(y,z+u)-f(y,z)}{u}uG(u)=0,
\end{split}
\eenn
since $f$ is differentiable and $\lim_{u\rightarrow 0}uG(u)=0$.
The last equality follows from the fact that if
\ebnn
\limsup_{u\rightarrow 0}uG(u)> 0
\eenn
then there exist a decreasing sequence $s_n\rightarrow 0$ and $d>0$ such that for each $n \in \N$
\ebnn
s_nG(s_n)>d.
\eenn
We can now find a subsequence $s_{n_k}$ satisfying
\ebnn
s_{n_k}\nu((s_{n_{k+1}},s_{n_k}])>\frac{d}{2}
\eenn
for each $k \in \N$. We reached a contradiction with $\int_{(0,\infty)}\min\{1,x\}\nu(dx)<\infty$. Thus the limit equals 0.
 Taking into account the above calculations, we follow the reasoning from \cite{MGZ}. That gives us for the last summand in  \eqref{adjoint_relation}
\eb
\label{adj_4_summand}
\begin{split}
\int_{\R^2} \int_0^\infty &\left[f(y,z+u) -f(y,z) \right]  p_t(y,z)\nu(du)dydz\\
&=-\int_{-\R^2}  f(y,w) \frac{\partial}{\partial w}\int_{0}^w G(w-z)p_t(y,z)dzdydw.
\end{split}
\ee
Let us introduce the following operator
\eb
\label{theta_operator}
\Theta_w g(w)=\int_0^w G(w-z)g(z)dz.
\ee
Laplace transform of its kernel equals
\eb
\label{G_laplace}
\lap[G](u)=\int_0^\infty e^{-ut}G(t)dt=\int_0^\infty \int_{(t,\infty)}e^{-ut}\nu(ds)dt=\int_0^\infty \frac{1}{u}(1-e^{-us})\nu(ds)=\frac{\Psi(u)}{u}.
\ee
Combining Eqs. \eqref{generator_jumps_end} and \eqref{adj_4_summand} we obtain the explicit formula for the adjoint operator $A^+$ with $p_t(y,z)$ in its domain,
and from Eq. \eqref{q_adjoint} we finally get
\eb
\begin{split}
\label{p_t_evolution}
\frac{\partial}{\partial t}p_t(y,z)=&- \frac{\partial}{\partial y}\left(F(y,z)p_t(y,z)\right)+\frac{\partial^2}{\partial y^2} \left(\frac{1}{2}\sigma^2(y,z)p_t(y,z)\right)\\
&+\int_{\R\setminus\{0\}}\left(p_t(y-E(z)x,z)
-p_t(y,z)
+E(z)x\frac{\partial}{\partial y}p_t(y,z)\ind_{\{\abs{x}<1\}}\right)\mu(dx)\\
&-\frac{\partial}{\partial z}\Theta_z p_t(y,z).
\end{split}
\ee
We derived the equation for temporal evolution of $p_t(y,z)$ which was our goal in the first step of the proof.
Let $w(x,t)$ denote the probability density function of the process $X(t)$.
In the next step we will find the relation between the densities $p_t(y,z)$ and $w(x,t)$. We start with denoting paths of the analyzed processes as $X(t,\omega)$ and $(Y(t,\omega), Z(t,\omega))$ for each $\omega \in \Omega$.
For each fixed interval $I$, similarly as in \cite{MGZ, Straka}, we define an auxiliary function
\eb
H_t(s,\omega,u)=\left\{\begin{array}{ll}
\ind_I(Y^-(s,\omega)) & \textrm{if $ Z^-(s,\omega)\leq t \leq Z^-(s,\omega)+u$} \\
0 & \textrm{otherwise} \end{array} \right .
\ee
and follow, with appropriate changes, the reasoning from the mentioned papers. This includes the observation
\eb\label{delta}
\ind_I(X(t,\omega))=\sum_{s>0}H_t(s,\omega,\Delta Z(s,\omega)),
\ee
where $\Delta Z(s)=Z(s)-Z^-(s)$. In our case the above equation is also valid since we excluded the case of compound Poisson processes and jumping times of $Z(t,\omega)$
are dense in $[0,\infty)$ almost surely. The second key equation is
 the compensation formula (Ch. XII, Proposition (1.10) in \cite{Revuz})
\eb
\label{compensation_formula}
\ex\left[\sum_{s>0}H_t(s,\omega,\Delta Z(s,\omega))\right]=\ex\left[\int_0^\infty \int_{0}^\infty H_t(s,\omega,u)\nu(du)ds \right].
\ee
  The difference between our case and \cite{MGZ} is that in the definition of $H_t$, instead of $Y(s,\omega)$, we put $Y^-(s,\omega)$ . Due to the L\'evy noise the previous process may no longer be continuous, but $Y^-$ is still left-continuous and therefore predictable, which is a condition for using the compensation formula. Another condition is that
  \ebnn
\sum_{s\in\R^+: \Delta Z(s)=0} H_t(s,\omega,\Delta Z(s,\omega))=0 \quad \mbox{almost surely},
\eenn
which is also fulfilled.
Indeed, assume to the contrary that for $\omega \in W \subset\Omega$ this sum does not vanish and $P(W)\neq 0$.
This means that there exist $s>0$ such that $Z^-(s,\omega)=t=Z(s,\omega)$.
However the paths of the process $Z(s,\omega)$ hit the single point $t$ with probability 0 (see \cite{Sato}),
which contradicts the assumption.
Consequently,
\eb
\label{w}
w(x,t)=\int_0^\infty \Theta_t p_{s}(x,t)ds.
\ee
Detailed derivation of the above result (with the Riemann-Liouville fractional integral instead of $\Theta_t$) is in \cite{MGZ}. Notice that
\eqb
\label{rozniczkowanie_theta}
\frac{d}{dt} \Theta_t f(t)= \Theta_t \frac{d}{dt}f(t)+f(0)G(t)
\eqe
for sufficiently smooth $f$. This can be proven calculating Laplace transform: for the left-hand side of
\eqref{rozniczkowanie_theta} we get
\eb
\nonumber
\mathcal{L}\left[\frac{d}{dt} \Theta_t f(t)\right](s)=s\mathcal{L}\left[ \Theta_t f(t)\right](s)=\Psi(s)\mathcal{L}[f(t)](s),
\ee
whereas the right-hand side equals
\eb
\nonumber
\begin{split}
&\mathcal{L}\left[\Theta_t \frac{d}{dt}f(t)+f(0)G(t)\right](s)=\frac{\Psi(s)}{s}\mathcal{L}\left[\frac{d}{dt}f(t)\right](s)+f(0)\frac{\Psi(s)}{s}=\Psi(s)\mathcal{L}[f(t)](s).
\end{split}
\ee
Thus
\eqb
\nonumber
\label{rozniczkowanie_theta_q}
\frac{\partial}{\partial t} \Theta_t p_s(x,t)= \Theta_t \frac{\partial}{\partial t}p_s(x,t)+p_s(x,0)G(t)
\eqe
and we have the following approximation
\eb
\label{szacowanie_rozbicie_na_2_calki}
\begin{split}
\int_0^{t_0}&\int_0^\infty\abs{\frac{\partial}{\partial t} \Theta_t p_s(x,t)}dsdt\\
&\leq\int_0^{t_0}\int_0^\infty\abs{\Theta_t \frac{\partial}{\partial t}p_s(x,t)}dsdt +\int_0^{t_0}\int_0^\infty p_s(x,0)G(t)dsdt.
\end{split}
\ee
We deal with both integrals separately. For the first one we observe that
\eb
\label{oszacowanie_calka_g}
\int_0^{t_0}G(u)du=\int_0^{t_0}\int_{(u,\infty)}\mu(dw)du=\int_{(0,\infty)}\min{(w,t_0)}\mu(dw)=K<\infty,
\ee
where $K>0$, and use Fubini's theorem together with the assumption  \eqref{dt_integrable}
\eb
\begin{split}
\int_0^{t_0}&\int_0^\infty\abs{\Theta_t\frac{\partial}{\partial t}  p_s(x,t)}dsdt\\
&\leq\int_0^{t_0}\int_0^\infty\int_0^t G(t-u) \abs{\frac{\partial}{\partial u} p_s(x,u)}dudsdt\\
&\leq K
 \int_0^{t_0}\int_0^\infty\abs{\frac{\partial}{\partial u} p_s(x,u)}dsdu
<\infty.
\end{split}
\ee
For the second integral in \eqref{szacowanie_rozbicie_na_2_calki} we obtain
\eb
\begin{split}
\int_0^{t_0}\int_0^\infty p_s(x,0)G(t)dsdt=K\int_0^\infty p_s(x,0)ds<\infty.
\end{split}
\ee
The last inequality is a consequence of Theorem 35.4 in \cite{Sato}.
Combining the approximations for both integrals we get
\eb
\label{gt_u_szac_rozniczkowanie2}
\int_0^{t_0}\int_0^\infty\abs{\frac{\partial}{\partial t} \Theta_t p_s(x,t)}dsdt<\infty.
\ee
Therefore after differentiating Eq. \eqref{w} with respect to $t$ we can move the derivative  inside the integral:
\eb
\label{w_diff}
\frac{\partial}{\partial t}w(x,t)=\int_0^\infty \frac{\partial}{\partial t} \Theta_t p_{s}(x,t)ds.
\ee
Now, applying  Eq. \eqref{p_t_evolution},
taking into account the facts that $\lim_{s\rightarrow \infty}p_{s}(x,t)=0$ and $p_{0}(x,t)=\delta_{(0,0)}(x,t)$ - the 2-dimensional Dirac delta, we obtain for $t>0$
\eb
\label{w_t_transformed}
\begin{split}
\frac{\partial}{\partial t}w(x,t)&=\int_0^\infty \left[\frac{\partial^2}{\partial x^2} \left(\frac{1}{2}\sigma^2(x,t)p_{s}(x,t)\right)- \frac{\partial}{\partial x}\left(F(x,t)p_{s}(x,t)\right)\right.\\
&\left.\quad+\int_{\R\setminus\{0\}}\left(p_s(x-E(t)y,t)
-p_s(x,t)
+E(t)y\frac{\partial}{\partial x}p_s(x,t)\ind_{\{\abs{y}<1\}}\right)\nu(dy)\right.\\
&\quad\left.-\frac{\partial}{\partial s}p_{s}(x,t)\right]ds \\
&=\frac{\partial^2}{\partial x^2}\left(\frac{1}{2}\sigma^2(x,t)\int_0^\infty  p_{s}(x,t)ds\right)- \frac{\partial}{\partial x}\left(F(x,t)\int_0^\infty p_{s}(x,t)ds\right)\\
&\quad+\int_{\R\setminus\{0\}}\left[\int_0^\infty p_s(x-E(t)y,t)ds
-\int_0^\infty p_s(x,t)ds\right.\\
&\left.\quad\quad+E(t)y\frac{\partial}{\partial x}\left( \int_0^\infty p_s(x,t)ds\right)\ind_{\{\abs{y}<1\}}\right]\nu(dy).
\end{split}
\ee
We can interchange the derivative and the integral here because of the assumptions \eqref{dx_integrable}. Changing the order of integration, $ds$ with $\nu(dy)$, is justified based on Fubini's theorem and the assumption \eqref{ds_integrable}.
Next we apply Fubini's theorem again, this time to Eq. \eqref{w}, obtaining
\eb
\begin{split}
w(x,t)=\int_0^\infty\int_0^t G(t-z)p_s(x,z)dzds=\Theta_t\int_0^\infty p_s(x,t)ds,
\end{split}
\ee
because the function $G(t-z)p_s(x,z)$ is non-negative. Therefore
\eb
\label{theta_inverse}
\int_0^\infty p_{s}(x,t)ds=\Theta_t^{-1}w(x,t),
\ee
and
\eqb
\label{theta_inverse2}
\int_0^\infty p_{s}(x-E(t)y,t)ds=(\Theta_t^{-1}w(r,t))\bigg|_{r=x-E(t)y},
\eqe
where $\Theta_t^{-1}$ is the left-inverse of the operator $\Theta_t$ defined in \eqref{theta_operator}.  It turns out that that $\Theta_t^{-1}=\Phi_t$ (see Eq. \eqref{operator}). Indeed, one can easily notice that $\Theta_t$ is a convolution operator with the kernel $G(u)$ and similarly, $\Phi_t$ is a composition of a derivative operator and a convolution operator with the kernel $M(u)$. Therefore, using Eqs. \eqref{G_laplace} and \eqref{kernel}  we get
\ebnn
\label{2_laplace}
\lap[\Phi_t\Theta_t f(t)](u)=u\lap[f](u)\lap[M](u)\lap[G](u)=u\lap[f](u)\frac{1}{\Psi(u)}\frac{\Psi (u)}{u}=\lap[f](u).
\eenn
 for a sufficiently smooth function $f$, which implies that $\Phi_t\Theta_t f(t)=f(t)$. Consequently, from Eqs. \eqref{w_t_transformed} and \eqref{theta_inverse} we obtain the desired result.
\end{proof}
The proof can not be easily extended to the case where the jump coefficient E is space dependent. In such situation the operator $A^+$ does not exist.
\begin{remark}
The above theorem can be used to approximate solutions of FFPE \eqref{ffp_general}
using Monte Carlo methods based on realizations of the process $X(t)$. Indeed, to simulate trajectories of $X(t)$, one only needs
to simulate the process $Y(t)$ (using the standard Euler scheme \cite{Janicki_Weron}) and the inverse subordinator
$S_\Psi(t)$ \cite{Magdziarz1,Magdziarz3}.

It also opens the possibility of analyzing fractional Cauchy problems \cite{Meer2,Meer2a,Meer9,Leonenko} in the general setting of ID subordinators and arbitrary space-time-dependent drift and diffusion coefficients.
\end{remark}

\section*{Acknowledgements}
This research was partially supported by NCN Maestro grant no. 2012/06/A/ST1/00258.

\bibliographystyle{amsplain}

\end{document}